\documentclass[11pt]{article}

\usepackage{amssymb}
\usepackage{amsmath}
\usepackage[utf8]{inputenc}
\usepackage{geometry}
\usepackage{graphicx}
\usepackage{amssymb}
\usepackage{listings}
\usepackage{epsfig}
\allowdisplaybreaks

\usepackage{float}
\usepackage{epstopdf}

\usepackage{tcolorbox}
\usepackage{tikz} 
\usepackage{pgfplots}
\usetikzlibrary{pgfplots.groupplots} 
\usetikzlibrary{arrows.meta, positioning, shadows}
\usepackage{circuitikz}
\usepackage{tikz}
\usepackage{circuitikz}
\usepackage{ragged2e} 
\justifying
\usetikzlibrary{matrix}
\usepackage{algorithm}
\usepackage[noend]{algpseudocode}

\usepackage{subcaption}

\title{The Edwards Model for fBm Loops and Starbursts}
\author{
	Wolfgang Bock and Yashika Jayathunga\\
	{\small Technomathematics Group}\\
	{\small	University of Kaiserslautern}\\
	{\small	P.\ O.\ Box 3049, 67653 Kaiserslautern, Germany}\\
	{\small E-Mail:bock@mathemaik.uni-kl.de, yashika.jayathunga@gmail.com}}

\title{An ADI Scheme for Two-sided Fractional Reaction-Diffusion Equations and Applications to an Epidemic Model}

\begin{document}        
        \maketitle
        \begin{abstract}
        Reaction-diffusion equations are often used in epidemiological models. In this paper we generalize the algorithm of Meerschaert and Tadjeran for fractional advection-dispersion flow equations to a coupled system of fractional reaction-diffusion like equation which arise from vector bourne disease modeling.
        \end{abstract}

\section{Introduction}
\label{sec:intro}
The modeling and understanding of infectious diseases is for many decades an object of intensive study. Going back to the classical SIR model from Kermack and McKendrick \cite{Kermack} which describes the time evolution of the number of susceptible (S), infected (I) and recovered (R) individuals by a system of ordinary differential equations various refinements were developed and extended exhaustively in the last 90 years. Among those extensions are the introduction of new compartments e.g.~to model vector-borne diseases such as Dengue or malaria, as well as more involved deterministic and stochastic models, see e.g.~\cite{PGSW17, GSW15, GRW15, GSW14, KAS14, RMT14, GNSW13, AKMS12, AGS12,ABC16, CG12, CEV09}.
Spatial disease spread can be modeled either in a discrete or continuous way. Popular space-discrete models are the metapopulation approach ~\cite{Arino,Dhirasakdanon,Jansen} and for Dengue recently \cite{BJ18, BJ19}, Cellular automata ~\cite{Yazhi,Slimi},  epidemic spatial networks~\cite{Bauch,BauchII,Newman,Colizza},  and lattice epidemic models ~\cite{Rhodes,Sazonov}.  
For space-continuous models integro-differential equation epidemic models ~\cite{Kot,Ralynn} and diffusion epidemic models ~\cite{Kunisch,Murray,Tuncer}, are studied.
In the last 10 years numerous fractional epidemic models~\cite{Brockmann, PINTO, Gonzalez, Nico9} were established.  A distinct feature of fractional derivatives is the capability to model long-range interactions. In a popular model, the second derivative in a classical diffusion model is substituted by $\alpha$- order derivative.
       
In this article, a fractional diffusion model is derived from the $SIRUV$ compartmental model with migration. To simulate this equation numerically  we adapt the  Alternating Directions Implicit (ADI) scheme with a Crank-Nicholson discretization to the fractional case.
For this purpose a shifted version of the typical Gr\"unwald-Letnikov finite difference approximation is used.  The ADI method joined together with a fractional Crank-Nicholson scheme for fractional diffusion examples was already implemented by Meerschaert et. al.~\cite{MR2284325,meerschaert2004finite,Meerschaert1,Meerschaert2}. The novelty of this article is that it is generalized to a system of coupled fractional reaction-diffusion equations. For this we derive the ADI splits with the corresponding Gr\"unwald-Letnikov operators.
A numerical scenario and a comparison with the classical diffusion case for Dirichlet boundary conditions can be found at the end of the article. 

\section{Model Definition}
In this article, the system of ordinary differential equations (ODEs)  for $SIV$ model is taken derive the fractional model. Instead of using the system of equations for $SIRUV$ model as in \cite{MR3805116}, a reduced form is used by using the simplification $R(t)=N-S(t)-I(t)$ and $U(t)=M-V(t)$ is given by the system of equations  (\ref{modelSIRUV}).  The corresponding system of ODEs is given as follows:
\begin{equation}
    \begin{aligned}
    \frac{dS(t)}{dt}  &= \mu \cdot (1-S(t))  -\beta \cdot S(t) \cdot V(t) =g_S \\
    \frac{dI(t)}{dt}  &=\beta \cdot S(t) \cdot V(t) -(\mu+\gamma)\cdot I(t)  =g_I\\
    \frac{dV(t)}{dt}  &=\vartheta \cdot (1-V(t)) \cdot I(t) -\nu \cdot  V(t)  =g_V\\
    \label{modelSIRUV}    
    \end{aligned}
    \end{equation}   
    where  $\beta$ and $\vartheta$ are the infection rate from vectors to hosts and hosts to vectors respectively.  The recovery rate from the compartment $I$  is given by $\gamma$. The birth and death rates of the hosts are equal and denoted by $\mu$ in order to maintain constant population size.  Similarly, a constant population of vectors is maintained by assuming that birth and death rates $\nu$ of the vectors to be equal.
    The initial conditions are given by   $S(0)$, $I(0)$ and $V(0)$ for the corresponding compartments.
    
    The susceptible and infected individuals are spatially distributed,  where $S(x,y,t)$, $I(x,y,t)$ and $V(x,y,t)$ represent the three state variables for the compartments. The initial conditions are given by the notations   $S(x,y,t_0)$, $I(x,y,t_0)$ and $V(x,y,t_0)$. The two-dimensional spatial variables are denoted by  $x$ and $y$. Model (\ref{modelSIRUV}) is redefined and can be written in a form of  a reaction-diffusion model as follows:
 \begin{align}
&    \dfrac{\partial S(x,y,t)}{\partial t} = g_S +a^S  \left \{     \dfrac{\partial^{2} S(x,y,t)}{\partial x^{2}} \right \}+b^S  \left \{  \dfrac{\partial^{2} S(x,y,t)}{\partial y^{2}} \right \}  \nonumber \\
&    \dfrac{\partial I(x,y,t)}{\partial t}= g_I +a^I \left \{  \dfrac{\partial^{2} I(x,y,t)}{\partial x^{2}} \right \}  + b^I \left \{  \dfrac{\partial^{2} I(x,y,t)}{\partial y^{2}} \right \} \nonumber  \\
&    \dfrac{\partial V(x,y,t)}{\partial t}=  g_V  +a^V\left \{  \dfrac{\partial^{2} V(x,y,t)}{\partial x^{2}} \right \}  + b^V \left \{  \dfrac{\partial^{2} V(x,y,t)}{\partial y^{2}} \right \}    \label{chapter6_eq35a}
\end{align}
 on a finite rectangular domain $x_L<x<x_H$ and $y_L<y<y_H$. The   fractional orders are given by  $1< \alpha_1 \leq 2$ and $1< \alpha _2 \leq 2$.   Dirichlet boundary conditions are used on the  boundary  $x_L\leq x \leq x_H$ and $y_L \leq y \leq y_H$:
    \begin{align*}
    & S(x=x_L,y,t)=
         S(x=x_R,y,t)=
     S(x,y=y_L,t)=
       S(x,y=y_R,t)= 0\nonumber \\
    &I(x=x_L,y,t)= 
        I(x=x_R,y,t)=  
         I(x,y=y_L,t)= 
         I(x,y=y_R,t)= 0\\
    &V(x=x_L,y,t)= 
        V(x=x_R,y,t)=  
         V(x,y=y_L,t)= 
         V(x,y=y_R,t)= 0.
    \end{align*}    \label{chapter6_eq35b}
    The fractional derivatives of the previous equations are replaced by two-sided fractional derivatives and hence, the two-sided fractional diffusion $SIV$-model yields,
     \begin{align}
 \dfrac{\partial S(x,y,t)}{\partial t}  &= g_S+a^S  \left \{     (1-r_1)\dfrac{\partial^{\alpha_1} S(x,y,t)}{\partial(-x)^{\alpha_1}}  +  r_1\dfrac{\partial^{\alpha_1} S(x,y,t)}{\partial x^{\alpha_1}} \right \} \nonumber\\
  &+ b^S  \left \{(1-r_2) \dfrac{\partial^{\alpha_2} S(x,y,t)}{\partial (-y)^{\alpha_2}}  +r_2 \dfrac{\partial^{\alpha_2} S(x,y,t)}{\partial y^{\alpha_2}} \right \}  
   \nonumber \\ 
 \dfrac{\partial I(x,y,t)}{\partial t}&=g_I+a^I \left \{     (1-r_1)\dfrac{\partial^{\alpha_1} I(x,y,t)}{\partial(-x)^{\alpha_1}}  +  r_1\dfrac{\partial^{\alpha_1} I(x,y,t)}{\partial x^{\alpha_1}} \right \} \nonumber\\ 
&+ b^I \left \{(1-r_2) \dfrac{\partial^{\alpha_2} I(x,y,t)}{\partial (-y)^{\alpha_2}}  +r_2 \dfrac{\partial^{\alpha_2} I(x,y,t)}{\partial y^{\alpha_2}} \right \}  \nonumber\\
\dfrac{\partial V(x,y,t)}{\partial t}&=g_V+a^V\left \{     (1-r_1)\dfrac{\partial^{\alpha_1} V(x,y,t)}{\partial(-x)^{\alpha_1}}  +  r_1\dfrac{\partial^{\alpha_1} V(x,y,t)}{\partial x^{\alpha_1}} \right \}\nonumber\\
&+ b^V \left \{(1-r_2) \dfrac{\partial^{\alpha_2} V(x,y,t)}{\partial (-y)^{\alpha_2}}  +r_2 \dfrac{\partial^{\alpha_2} V(x,y,t)}{\partial y^{\alpha_2}} \right \} \label{chapter6_eq35C},
\end{align} 
with weights $r_1, r_2 \in [0,1]$, where $\dfrac{\partial^{\alpha_i} F(x,y,t)}{\partial(-x)^{\alpha_i}}$ and $\dfrac{\partial^{\alpha_i} F(x,y,t)}{\partial(-y)^{\alpha_i}}$ denote the negative (right) fractional derivatives. 
    
\section{Numerical Scheme}
 A Crank-Nicholson type system of finite difference equations can be obtained by substituting the shifted Gr{\"u}nwald into the differential equation centered at time $t_{n+1/2}=\dfrac{1}{2} (t_{n+1}+t_n)$.
    
\begin{align}
S_{i,j}^{n+1}-S_{i,j}^{n}  &= -\Delta  t\left \{\mu-\mu S_{i,j}^{n+1/2}-\beta S_{i,j}^{n+1/2} V_{i,j}^{n+1/2}  \right\} \nonumber\\
&\quad +\left. {} \dfrac{\Delta  t}{2} \left \{ (1-r_1) (\delta_{\alpha_1 x}^{S^-} S_{i,j}^{n+1} + \delta_{\alpha_1 x}^{S^-} S_{i,j}^{n} ) + r_1 (\delta_{\alpha_1 x}^{S^+} S_{i,j}^{n+1} + \delta_{\alpha_1 x}^{S^+} S_{i,j}^{n} ) \right.\right\}  \nonumber \\ 
&\quad+ \dfrac{\Delta  t}{2} \left \{(1-r_2) (\delta_{\alpha_2 y}^{S^-} S_{i,j}^{n+1} + \delta_{\alpha_2 y}^{S^-} S_{i,j}^{n} ) + r_2 (\delta_{\alpha_2 y}^{S^+} S_{i,j}^{n+1} + \delta_{\alpha_2 y}^{S^+} S_{i,j}^{n} )\right\} \nonumber\\
I_{i,j}^{n+1}-I_{i,j}^{n}  &= \Delta  t \left \{\beta S_{i,j}^{n+1/2} V_{i,j}^{n+1/2}-(\mu+\gamma) I_{i,j}^{n+1/2}  \right\} \nonumber\\
&\quad + \dfrac{\Delta  t}{2} \left\{ (1-r_1) (\delta_{\alpha_1 x}^{I^-} I_{i,j}^{n+1} + \delta_{\alpha_1 x}^{I^-} I_{i,j}^{n} ) + r_1 (\delta_{\alpha_1 x}^{I^+} I_{i,j}^{n+1} + \delta_{\alpha_1 x}^{I^+} I_{i,j}^{n} ) \right\}  \nonumber\\ 
&\quad+ \dfrac{\Delta  t}{2} \left \{(1-r_2) (\delta_{\alpha_2 y}^{I^-} I_{i,j}^{n+1} + \delta_{\alpha_2 y}^{I^-} I_{i,j}^{n} ) + r_2 (\delta_{\alpha_2 y}^{I^+} I_{i,j}^{n+1} + \delta_{\alpha_2 y}^{I^+} I_{i,j}^{n} )\right\}\\
V_{i,j}^{n+1}-V_{i,j}^{n}  &= \Delta  t \left \{\vartheta  I_{i,j}^{n+1/2}-\vartheta  V_{i,j}^{n+1/2} I_{i,j}^{n+1/2}-\nu V_{i,j}^{n+1/2}  \right\} \nonumber\\
&\quad + \dfrac{\Delta  t}{2} \left \{ (1-r_1) (\delta_{\alpha_1 x}^{V^-} V_{i,j}^{n+1} + \delta_{\alpha_1 x}^{V^-} V_{i,j}^{n} ) + r_1 (\delta_{\alpha_1 x}^{V^+} V_{i,j}^{n+1} + \delta_{\alpha_1 x}^{V^+} V_{i,j}^{n} ) \right\}  \nonumber \\ 
&\quad+ \dfrac{\Delta  t}{2} \left \{(1-r_2) (\delta_{\alpha_2 y}^{V^-} V_{i,j}^{n+1} + \delta_{\alpha_2 y}^{V^-} V_{i,j}^{n} ) + r_2 (\delta_{\alpha_2 y}^{V^+} V_{i,j}^{n+1} + \delta_{\alpha_2 y}^{V^+} V_{i,j}^{n} )\right\}\nonumber 
 \end{align}
 
After rearranging the terms, the previous equation can be written in the operator notations as (\ref{eq5x}):
\begin{multline}
(1-\dfrac{\Delta  t}{2} \{(1-r_1)\delta_{\alpha_1 x}^{S^-} +r_1\delta_{\alpha_1 x}^{S^+}\}-\dfrac{\Delta  t}{2} \{(1-r_2)\delta_{\alpha_2 y}^{S^-} +r_2\delta_{\alpha_2 y}^{S^+}\})S_{i,j}^{n+1}  \nonumber\\
    =(1+\dfrac{\Delta  t}{2} \{(1-r_1)\delta_{\alpha_1 x}^{S^-} +r_1\delta_{\alpha_1 x}^{S^+}\}+\dfrac{\Delta  t}{2} \{(1-r_2)\delta_{\alpha_2 y}^{S^-} +r_2\delta_{\alpha_2 y}^{S^+}\})S_{i,j}^{n}\\ 
-\Delta  t \left \{\mu-\mu S_{i,j}^{n+1/2}-\beta S_{i,j}^{n+1/2} V_{i,j}^{n+1/2}  \right\}
\end{multline}
\begin{multline}
(1-\dfrac{\Delta  t}{2} \{(1-r_1)\delta_{\alpha_1 x}^{I^-} +r_1\delta_{\alpha_1 x}^{I^+}\}-\dfrac{\Delta  t}{2} \{(1-r_2)\delta_{\alpha_2 y}^{I^-} +r_2\delta_{\alpha_2 y}^{I^+}\})I_{i,j}^{n+1} \nonumber \\
    =(1+\dfrac{\Delta  t}{2} \{(1-r_1)\delta_{\alpha_1 x}^{I^-} +r_1\delta_{\alpha_1 x}^{I^+}\}+\dfrac{\Delta  t}{2} \{(1-r_2)\delta_{\alpha_2 y}^{I^-} +r_2\delta_{\alpha_2 y}^{I^+}\})I_{i,j}^{n}\\ 
+\Delta  t \left \{\beta S_{i,j}^{n+1/2} V_{i,j}^{n+1/2}-(\mu+\gamma) I_{i,j}^{n+1/2}  \right\}
\end{multline}
\begin{multline}
(1-\dfrac{\Delta  t}{2} \{(1-r_1)\delta_{\alpha_1 x}^{V^-} +r_1\delta_{\alpha_1 x}^{V^+}\}-\dfrac{\Delta  t}{2} \{(1-r_2)\delta_{\alpha_2 y}^{V^-} +r_2\delta_{\alpha_2 y}^{V^+}\})V_{i,j}^{n+1}  \\
    =(1+\dfrac{\Delta  t}{2} \{(1-r_1)\delta_{\alpha_1 x}^{V^-} +r_1\delta_{\alpha_1 x}^{V^+}\}+\dfrac{\Delta  t}{2} \{(1-r_2)\delta_{\alpha_2 y}^{V^-} +r_2\delta_{\alpha_2 y}^{V^+}\})V_{i,j}^{n}\\ 
+\Delta  t \left \{\vartheta  I_{i,j}^{n+1/2}-\vartheta  V_{i,j}^{n+1/2} I_{i,j}^{n+1/2}-\nu V_{i,j}^{n+1/2}  \right\}\label{eq5x}
\end{multline}
 Multi-dimensional diffusion equations are often solved with alternating directions implicit methods (ADI), where splitting is used to significantly reduce the computational work \cite{ADI}. These techniques use a perturbation of Equation (\ref{eq5x}) in order to derive schemes that requires only the implicit numerical solution in one direction where the other spatial direction is computed iteratively. We obtain the equations\\
    for $S$,
     \begin{equation}
    \begin{split}
    (1-\dfrac{\Delta  t}{2}\{ (1-r_1)\delta_{\alpha_1 x}^{S^-} +r_1\delta_{\alpha_1 x}^{S^+}\})(1-\dfrac{\Delta  t}{2}\{(1-r_2)\delta_{\alpha_2 y}^{S^-} +r_2\delta_{\alpha_2 y}^{S^+}\}) S_{i,j}^{n+1}\\
    =(1+\dfrac{\Delta  t}{2}\{(1-r_1)\delta_{\alpha_1 x}^{S^-} +r_1\delta_{\alpha_1 x}^{S^+}\}(1+\dfrac{\Delta  t}{2} \{(1-r_2)\delta_{\alpha_2 y}^{S^-} +r_2\delta_{\alpha_2 y}^{S^+}\}) S_{i,j}^{n} \\-\Delta  t \left \{\mu-\mu S_{i,j}^{n+1/2}-\beta S_{i,j}^{n+1/2} V_{i,j}^{n+1/2}  \right\},
    \end{split} \label{chapter6_eq20}
    \end{equation}
    for $I$,
    \begin{equation}
    \begin{split}
    (1-\dfrac{\Delta  t}{2}\{ (1-r_1)\delta_{\alpha_1 x}^{I^-} +r_1\delta_{\alpha_1 x}^{I^+}\})(1-\dfrac{\Delta  t}{2}\{ (1-r_2)\delta_{\alpha_2 y}^{I^-} +r_2\delta_{\alpha_2 y}^{I^+}\}) I_{i,j}^{n+1}\\
    =(1+\dfrac{\Delta  t}{2}\{ (1-r_1)\delta_{\alpha_1 x}^{I^-} +r_1\delta_{\alpha_1 x}^{I^+}\})(1+\dfrac{\Delta  t}{2}\{ (1-r_2)\delta_{\alpha_2 y}^{I^-} +r_2\delta_{\alpha_2 y}^{I^+}\}) I_{i,j}^{n} \\+\Delta  t \left \{\beta S_{i,j}^{n+1/2} V_{i,j}^{n+1/2}-(\mu+\gamma) I_{i,j}^{n+1/2}  \right\},
    \end{split} \label{chapter6_eq21}
    \end{equation}
    for $V$,
    \begin{equation}
    \begin{split}
    (1-\dfrac{\Delta  t}{2}\{ (1-r_1)\delta_{\alpha_1 x}^{V^-} +r_1\delta_{\alpha_1 x}^{V^+}\})(1-\dfrac{\Delta  t}{2}\{ (1-r_2)\delta_{\alpha_2 y}^{V^-} +r_2\delta_{\alpha_2 y}^{I^+}\}) V_{i,j}^{n+1}\\
    =(1+\dfrac{\Delta  t}{2}\{ (1-r_1)\delta_{\alpha_1 x}^{V^-} +r_1\delta_{\alpha_1 x}^{V^+}\})(1+\dfrac{\Delta  t}{2}\{ (1-r_2)\delta_{\alpha_2 y}^{V^-} +r_2\delta_{\alpha_2 y}^{I^+}\}) V_{i,j}^{n} \\+\Delta  t \left \{\vartheta  I_{i,j}^{n+1/2}-\vartheta  V_{i,j}^{n+1/2} I_{i,j}^{n+1/2}-\nu V_{i,j}^{n+1/2}.  \right\}.
    \end{split} \label{chapter6_eq21A}
    \end{equation}
 The equations (\ref{chapter6_eq20}), (\ref{chapter6_eq21}) and (\ref{chapter6_eq21A}) form Peaceman-Rachford type matrix equations defining ADI method.
 This can be  split as follows: \\
    for $S$,
    \begin{equation}
    \begin{split}
    (1-\dfrac{\Delta  t}{2}\{ (1-r_1)\delta_{\alpha_1 x}^{S^-} +r_1\delta_{\alpha_1 x}^{S^+}\}) S_{i,j}^{*}=(1+\dfrac{\Delta  t}{2}\{ (1-r_2)\delta_{\alpha_2 y}^{S^-} +r_2\delta_{\alpha_2 y}^{S^+}\})S_{i,j}^{n}\\ -\dfrac{\Delta  t}{2} \left \{\mu-\mu S_{i,j}^{n+1/2}-\beta S_{i,j}^{n+1/2} V_{i,j}^{n+1/2}  \right\}
    \end{split} \label{chapter6_eq22}
    \end{equation}
    \begin{equation}
    \begin{split}
    (1-\dfrac{\Delta  t}{2}\{ (1-r_2)\delta_{\alpha_2 y}^{S^-} +r_2\delta_{\alpha_2 y}^{S^+}\}) S_{i,j}^{n+1}=(1+\dfrac{\Delta  t}{2}\{ (1-r_1)\delta_{\alpha_1 x}^{S^-} +r_1\delta_{\alpha_1 x}^{S^+}\}) S_{i,j}^{*}\\-\dfrac{\Delta  t}{2}\left \{\mu-\mu S_{i,j}^{n+1/2}-\beta S_{i,j}^{n+1/2} V_{i,j}^{n+1/2}  \right\}
    \end{split} \label{chapter6_eq23}
    \end{equation}
    for $I$,
    \begin{equation}
    \begin{split}
    (1-\dfrac{\Delta  t}{2}\{ (1-r_1)\delta_{\alpha_1 x}^{I^-} +r_1\delta_{\alpha_1 x}^{I^+}\}) I_{i,j}^{*}=(1+\dfrac{\Delta  t}{2}\{ (1-r_2)\delta_{\alpha_2 y}^{I^-} +r_2\delta_{\alpha_2 y}^{I^+}\})I_{i,j}^{n} \\ +\dfrac{\Delta  t}{2} \left \{\beta S_{i,j}^{n+1/2} V_{i,j}^{n+1/2}-(\mu+\gamma) I_{i,j}^{n+1/2}  \right\}
    \end{split} \label{chapter6_eq24}
    \end{equation}
    \begin{equation}
    \begin{split}
    (1-\dfrac{\Delta  t}{2}\{ (1-r_2)\delta_{\alpha_2 y}^{I^-} +r_2\delta_{\alpha_2 y}^{I^+}\}) I_{i,j}^{n+1}=(1+\dfrac{\Delta  t}{2}\{ (1-r_1)\delta_{\alpha_1 x}^{I^-} +r_1\delta_{\alpha_1 x}^{I^+}\}) I_{i,j}^{*}\\ +\dfrac{\Delta  t}{2} \left \{\beta S_{i,j}^{n+1/2} V_{i,j}^{n+1/2}-(\mu+\gamma) I_{i,j}^{n+1/2}  \right\}.
    \end{split} \label{chapter6_eq25}
    \end{equation}
    for $V$,
    \begin{equation}
    \begin{split}
    (1-\dfrac{\Delta  t}{2}\{ (1-r_1)\delta_{\alpha_1 x}^{V^-} +r_1\delta_{\alpha_1 x}^{V^+}\}) V_{i,j}^{*}=(1+\dfrac{\Delta  t}{2}\{ (1-r_2)\delta_{\alpha_2 y}^{V^-} +r_2\delta_{\alpha_2 y}^{I^+}\})V_{i,j}^{n} \\ +\dfrac{\Delta  t}{2} \left \{\beta S_{i,j}^{n+1/2} V_{i,j}^{n+1/2}-(\mu+\gamma) I_{i,j}^{n+1/2}  \right\}
    \end{split} \label{chapter6_eq26a}
    \end{equation}
    \begin{equation}
    \begin{split}
    (1-\dfrac{\Delta  t}{2}\{ (1-r_2)\delta_{\alpha_2 y}^{I^-} +r_2\delta_{\alpha_2 y}^{I^+}\}) I_{i,j}^{n+1}=(1+\dfrac{\Delta  t}{2}\{ (1-r_1)\delta_{\alpha_1 x}^{I^-} +r_1\delta_{\alpha_1 x}^{I^+}\}) I_{i,j}^{*}\\ +\dfrac{\Delta  t}{2} \left \{\beta S_{i,j}^{n+1/2} V_{i,j}^{n+1/2}-(\mu+\gamma) I_{i,j}^{n+1/2}  \right\}.
    \end{split} \label{chapter6_eq27a}
    \end{equation}
    To observe this multiply (\ref{chapter6_eq22}), (\ref{chapter6_eq24}) and (\ref{chapter6_eq26a}) by $(1+\dfrac{\Delta  t}{2}( (1-r_1)\delta_{\alpha_1 x}^{-} +r_1\delta_{\alpha_1 x}^{+}))$ on both sides and then multiply (\ref{chapter6_eq23}), (\ref{chapter6_eq25}) and (\ref{chapter6_eq27a}) by $(1-\dfrac{\Delta  t}{2}( (1-r_1)\delta_{\alpha_1 x}^{-} +r_1\delta_{\alpha_1 x}^{+}))$ on both sides to obtain the Equations (\ref{chapter6_eq20}), (\ref{chapter6_eq21}) and (\ref{chapter6_eq21A}). Equations(\ref{chapter6_eq22}), (\ref{chapter6_eq23}), (\ref{chapter6_eq24}), (\ref{chapter6_eq25}), (\ref{chapter6_eq26a})   and (\ref{chapter6_eq27a})  calculates intermediate solutions $S_{i,j}^{*}$, $I_{i,j}^{*}$ and $V_{i,j}^{*}$ in order to develop the numerical solutions to $S$, $I$ and $V$ at time step $n$ to the numerical solution  $S_{i,j}^{n+1}$, $I_{i,j}^{n+1}$ and $V_{i,j}^{n+1}$ at time $t_{n+1}$. 
  
  The algorithm of ADI splitting method \cite{MR2284325} to solve the $SIV$-fractional diffusion model is given by:
   \begin{algorithm}[H]
      \caption{ADI Scheme}
      \begin{algorithmic}[1]
          \State     In order to acquire the intermediate solution slice  $S_{i,j}^*$, $I_{i,j}^*$ and $V_{i,j}^*$, a set of $N_x-1$ equations at the points $x_i$, where, $i=1,2,...,N_x-1$ defined by equation (\ref{chapter6_eq22}) and (\ref{chapter6_eq24}) are solved initially on each fixed horizontal slice $y=y_k$ where  $k=1,2,...,N_y-1$.    
          \State   In addition, by alternating  the spatial direction on every fixed verticle slice $x=x_k$ ($k=1,2,...,N_x-1$) a set of $N_y-1$ equations are solved at the points $y_j$ where $j=1,2,...,N_y-1$ defined by the equations (\ref{chapter6_eq23}) and (\ref{chapter6_eq25}) in order to obtain the solution for $S_{k,j}^{n+1}$, $I_{k,j}^{n+1}$ and $V_{k,j}^{n+1}$ at time $n+1$.          
      \end{algorithmic}    
  \end{algorithm}
    The shifted Gr{\"u}nwald operators used in this model yield
               \begin{align}
               &    \delta^{S^-}_{\alpha_1,x} S_{i,j}^n= \dfrac{(a_{ij}^S)}{(\Delta x)^{\alpha_1}} \sum_{k=0}^{N_x-i+1} g_{\alpha_1,k} \cdot  S_{i+k-1,j}^n
               &&    \delta^{S^+}_{\alpha_1, x} S_{i,j}^n= \dfrac{(a_{ij}^S)_n}{(\Delta x)^{\alpha_1}} \sum_{k=0}^{i+1} g_{\alpha_1,k} \cdot  S_{i-k+1,j}^n \nonumber \\
               &    \delta^{S^-}_{\alpha_2,y} S_{i,j}^n= \dfrac{(b_{ij}^S)}{(\Delta y)^{\alpha_2}} \sum_{k=0}^{N_y-j+1} g_{\alpha_2,k} \cdot  S_{i,j+k-1}^n
               &&    \delta^{S+}_{\alpha_2, y} S_{i,j}^n= \dfrac{(b_{ij}^S)}{(\Delta y)^{\alpha_2}} \sum_{k=0}^{j+1} g_{\alpha_2,k} \cdot  S_{i,j-k+1}^n \nonumber \\
               &    \delta^{I^-}_{\alpha_1,x} I_{i,j}^n= \dfrac{(a_{ij}^I)}{(\Delta x)^{\alpha_1}} \sum_{k=0}^{N_x-i+1} g_{\alpha_1,k} \cdot  I_{i+k-1,j}^n
               &&    \delta^{I^+}_{\alpha_1, x} I_{i,j}^n= \dfrac{(a_{ij}^I)_n}{(\Delta x)^{\alpha_1}} \sum_{k=0}^{i+1} g_{\alpha_1,k} \cdot  I_{i-k+1,j}^n \nonumber \\
               &    \delta^{I^-}_{\alpha_2,y} I_{i,j}^n= \dfrac{(b_{ij}^I)}{(\Delta y)^{\alpha_2}} \sum_{k=0}^{N_y-j+1} g_{\alpha_2,k} \cdot  I_{i,j+k-1}^n
               &&    \delta^{I+}_{\alpha_2, y} I_{i,j}^n= \dfrac{(b_{ij}^I)}{(\Delta y)^{\alpha_2}} \sum_{k=0}^{j+1} g_{\alpha_2,k} \cdot  I_{i,j-k+1}^n\nonumber \\
               &    \delta^{V^-}_{\alpha_1,x} V_{i,j}^n= \dfrac{(a_{ij}^V)}{(\Delta x)^{\alpha_1}} \sum_{k=0}^{N_x-i+1} g_{\alpha_1,k} \cdot  V_{i+k-1,j}^n
               &&    \delta^{V^+}_{\alpha_1, x} V_{i,j}^n= \dfrac{(a_{ij}^V)_n}{(\Delta x)^{\alpha_1}} \sum_{k=0}^{i+1} g_{\alpha_1,k} \cdot  V_{i-k+1,j}^n \nonumber \\
               &    \delta^{V^-}_{\alpha_2,y} V_{i,j}^n= \dfrac{(b_{ij}^V)}{(\Delta y)^{\alpha_2}} \sum_{k=0}^{N_y-j+1} g_{\alpha_2,k} \cdot  V_{i,j+k-1}^n
               &&    \delta^{V+}_{\alpha_2, y} V_{i,j}^n= \dfrac{(b_{ij}^V)}{(\Delta y)^{\alpha_2}} \sum_{k=0}^{j+1} g_{\alpha_2,k} \cdot  V_{i,j-k+1}^n 
               \end{align}    \label{chapter6_eq26}
                
Analogously finite difference schemes for the $S, I$ and $V $ compartments are obtained by substituting the shifted Gr{\"u}nwald operator into the equations from before. The corresponding ADI scheme reads:
\begin{itemize}
\item [](i). ADI split I: 
\begin{multline*}
S_{i,j}^{*}-r_1 D_{ij}^{S}  \sum_{k=0}^{i+1} g_{\alpha_1,k}\cdot  S_{i-k+1,j}^{*} 
-(1-r_1) D_{ij}^{S} \sum_{k=0}^{N_x-i+1} g_{\alpha_1,k}\cdot  S_{i+k-1,j}^{*} \\   
= S_{i,j}^{n}+r_2 E_{ij}^{S}  \sum_{k=0}^{j+1} g_{\alpha_2,k}\cdot  S_{i,j-k+1}^{n}
+(1-r_2) E_{ij}^{S}  \sum_{k=0}^{N_y-j+1} g_{\alpha_2,k}\cdot  S_{i,j+k-1}^{n} \nonumber \\ 
+ -\dfrac{\Delta  t}{2} \left \{  \mu-\beta S_{i,j}^{n+1/2} V_{i,j}^{n+1/2}-\mu S_{i,j}^{n+1/2}\right \} \nonumber, \\[.3cm]
\end{multline*}
\begin{multline*}
I_{i,j}^{*}-r_1  D_{ij}^{I} \sum_{k=0}^{i+1} g_{\alpha_1,k}\cdot  I_{i-k+1,j}^{*} 
- (1-r_1)D_{ij}^{I}  \sum_{k=0}^{N_x-i+1} g_{\alpha_1,k}\cdot  I_{i+k-1,j}^{*} \\  = I_{i,j}^{n}+r_2 E_{ij}^{I}  \sum_{k=0}^{j+1} g_{\alpha_2,k}\cdot  I_{i,j-k+1}^{n} 
+(1-r_2) E_{ij}^{I}  \sum_{k=0}^{N_y-j+1} g_{\alpha_2,k}\cdot  I_{i,j+k-1}^{n}  \\ 
+\dfrac{\Delta  t}{2} \left\{\beta S_{i,j}^{n+1/2} V_{i,j}^{n+1/2} -(\mu+\gamma)I_{i,j}^{n+1/2} \right \}, \\[.3cm] \end{multline*}
\begin{multline}
V_{i,j}^{*}-r_1  D_{ij}^{V} \sum_{k=0}^{i+1} g_{\alpha_1,k}\cdot  V_{i-k+1,j}^{*} 
- (1-r_1)D_{ij}^{V}  \sum_{k=0}^{N_x-i+1} g_{\alpha_1,k}\cdot  V_{i+k-1,j}^{*}\\   
= V_{i,j}^{n}+r_2 E_{ij}^{V}  \sum_{k=0}^{j+1} g_{\alpha_2,k}\cdot  V_{i,j-k+1}^{n} 
+(1-r_2) E_{ij}^{V}  \sum_{k=0}^{N_y-j+1} g_{\alpha_2,k}\cdot  V_{i,j+k-1}^{n} \nonumber \\ 
+\dfrac{\Delta  t}{2} \left\{ \vartheta (1-V_{i,j}^{n+1/2})  I_{i,j}^{n+1/2} -\nu V_{i,j}^{n+1/2}   \right \}, \\ 
\end{multline}
where $D_{ij}^{X}=\dfrac{\Delta t a_{ij}^X}{2 (\Delta x)^{\alpha_1}}$,    and  $E_{ij}^{X}=\dfrac{\Delta t b_{ij}^X}{2 (\Delta y)^{\alpha_2}}$ . X represents the compartments $S$, $I$ and $V$. 
    
\item [](ii). ADI split II:
\begin{multline*}
S_{i,j}^{n+1}-r_2 E_{ij}^{S}  \sum_{k=0}^{j+1} g_{\alpha_2,k}\cdot  S_{i,j-k+1}^{n+1}-(1-r_2) E_{ij}^{S} \sum_{k=0}^{N_y-j+1} g_{\alpha_2,k}\cdot  S_{i,j+k-1}^{n+1}\\   
= S_{i,j}^{*}+r_1 D_{ij}^{S}  \sum_{k=0}^{i+1} g_{\alpha_1,k}\cdot  S_{i-k+1,j}^{*} +(1-r_1) D_{ij}^{S}  \sum_{k=0}^{N_x-i+1} g_{\alpha_1,k}\cdot  S_{i+k-1,j}^{*} \\ 
-\dfrac{\Delta  t}{2} \left \{  \mu-\beta S_{i,j}^{n+1/2} V_{i,j}^{n+1/2}-\mu S_{i,j}^{n+1/2}\right \} \\[.2 cm]
\end{multline*}
\begin{multline*} 
I_{i,j}^{n+1}-r_2 E_{ij}^{I}  \sum_{k=0}^{j+1} g_{\alpha_2,k}\cdot  I_{i,j-k+1}^{n+1}-(1-r_2) E_{ij}^{I} \sum_{k=0}^{N_y-j+1} g_{\alpha_2,k}\cdot  I_{i,j+k-1}^{n+1}  \\
= I_{i,j}^{*}+r_1 D_{ij}^{I}  \sum_{k=0}^{i+1} g_{\alpha_1,k}\cdot  I_{i-k+1,j}^{*} +(1-r_1) D_{ij}^{I}  \sum_{k=0}^{N_x-i+1} g_{\alpha_1,k}\cdot  I_{i+k-1,j}^{*} \\ 
-\dfrac{\Delta  t}{2} \left\{\beta S_{i,j}^{n+1/2} V_{i,j}^{n+1/2} -(\mu+\gamma)I_{i,j}^{n+1/2} \right \} \\[.2 cm]
\end{multline*}
\begin{multline*} 
V_{i,j}^{n+1}-r_2 E_{ij}^{V}  \sum_{k=0}^{j+1} g_{\alpha_2,k}\cdot  V_{i,j-k+1}^{n+1}-(1-r_2) E_{ij}^{V} \sum_{k=0}^{N_y-j+1} g_{\alpha_2,k}\cdot  V_{i,j+k-1}^{n+1}   \\
= V_{i,j}^{*}+r_1 D_{ij}^{V}  \sum_{k=0}^{i+1} g_{\alpha_1,k}\cdot  V_{i-k+1,j}^{*}+(1-r_1) D_{ij}^{V}  \sum_{k=0}^{N_x-i+1} g_{\alpha_1,k}\cdot  V_{i+k-1,j}^{*} \\ 
-\dfrac{\Delta  t}{2} \left\{ \vartheta (1-V_{i,j}^{n+1/2})  I_{i,j}^{n+1/2} -\nu V_{i,j}^{n+1/2}   \right \} 
 \end{multline*}
\end{itemize}    
              Before solving the system of equations defined by ADI split I  and  ADI split II, the intermediate solutions $S_{ij}^*$, $I_{ij}^*$ and $V_{ij}^*$ must be treated with care on the boundary in order to  preserve the consistency of the set of equations defined by (\ref{chapter6_eq22}), (\ref{chapter6_eq23}), (\ref{chapter6_eq24}), (\ref{chapter6_eq25}), (\ref{chapter6_eq26a})   and (\ref{chapter6_eq27a})
              with (\ref{chapter6_eq20}), (\ref{chapter6_eq21}) and (\ref{chapter6_eq21A}).
              By subtracting (\ref{chapter6_eq23}) from (\ref{chapter6_eq22}),  (\ref{chapter6_eq25}) from (\ref{chapter6_eq24}) and (\ref{chapter6_eq27a}) from (\ref{chapter6_eq26a}) we obtain,
              
               \begin{equation}
               (1-\dfrac{\Delta  t}{2}\{ (1-r_2)\delta_{\alpha_2 y}^{S^-} +r_2\delta_{\alpha_2 y}^{S^+}\}) S_{i,j}^{n+1}+(1+\dfrac{\Delta  t}{2}\{ (1-r_2)\delta_{\alpha_2 y}^{S^-} +r_2\delta_{\alpha_2 y}^{S^+}\}) S_{i,j}^{n}=2 S_
               {i,j}^* \label{chapter6_eq29}
               \end{equation}
               \begin{equation}
               (1-\dfrac{\Delta  t}{2}\{ (1-r_2)\delta_{\alpha_2 y}^{I^-} +r_2\delta_{\alpha_2 y}^{I^+}\}) I_{i,j}^{n+1}+(1+\dfrac{\Delta  t}{2}\{ (1-r_2)\delta_{\alpha_2 y}^{I^-} +r_2\delta_{\alpha_2 y}^{I^+}\}) I_{i,j}^{n}=2 I_
               {i,j}^*
               \label{chapter6_eq30}
               \end{equation}
                   \begin{equation}
               (1-\dfrac{\Delta  t}{2}\{ (1-r_2)\delta_{\alpha_2 y}^{V^-} +r_2\delta_{\alpha_2 y}^{V^+}\}) V_{i,j}^{n+1}+(1+\dfrac{\Delta  t}{2}\{ (1-r_2)\delta_{\alpha_2 y}^{V^-} +r_2\delta_{\alpha_2 y}^{V^+}\}) V_{i,j}^{n}=2 V_
               {i,j}^*
 \label{chapter6_eq30a}
               \end{equation}
    The boundary conditions for the intermediate solutions $S_{ij}^*$, $I_{ij}^*$ and $V_{ij}^*$ ( $i.e.\ $, $i=0 \quad \text{or}\quad i=N_x \quad\text{for} \quad j=1, ..., N_y-1$ ) required to solve the set of equations  (\ref{chapter6_eq20}), (\ref{chapter6_eq21}) and (\ref{chapter6_eq21A}) are of the form
    \begin{align}
    2S_
    {0,j}^*=(1-\dfrac{\Delta  t}{2}\{ (1-r_2)\delta_{\alpha_2 y}^{S^-} +r_2\delta_{\alpha_2 y}^{S^+}\}) S_{0,j}^{n+1}+(1+\dfrac{\Delta  t}{2}\{ (1-r_2)\delta_{\alpha_2 y}^{S^-} +r_2\delta_{\alpha_2 y}^{S^+}\}) S_{0,j}^{n} \nonumber \\
   2 S_
    {N_x,j}^*=(1-\dfrac{\Delta  t}{2}\{ (1-r_2)\delta_{\alpha_2 y}^{S^-} +r_2\delta_{\alpha_2 y}^{S^+}\}) S_{N_x,j}^{n+1}+(1+\dfrac{\Delta  t}{2}\{ (1-r_2)\delta_{\alpha_2 y}^{S^-} +r_2\delta_{\alpha_2 y}^{S^+}\}) S_{N_x,j}^{n}
    \label{chapter6_eq31}
    \end{align}
    \begin{align}
    2I_
    {0,j}^*=(1-\dfrac{\Delta  t}{2}\{ (1-r_2)\delta_{\alpha_2 y}^{I^-} +r_2\delta_{\alpha_2 y}^{I^+}\}) I_{0,j}^{n+1}+(1+\dfrac{\Delta  t}{2}\{ (1-r_2)\delta_{\alpha_2 y}^{I^-} +r_2\delta_{\alpha_2 y}^{I^+}\}) I_{0,j}^{n}\nonumber \\
    2I_
    {N_x,j}^*=(1-\dfrac{\Delta  t}{2}\{ (1-r_2)\delta_{\alpha_2 y}^{I^-} +r_2\delta_{\alpha_2 y}^{I^+}\}) I_{N_x,j}^{n+1}+(1+\dfrac{\Delta  t}{2}\{ (1-r_2)\delta_{\alpha_2 y}^{I^-} +r_2\delta_{\alpha_{2} y}^{I^+}\}) I_{N_x,j}^{n}.
    \label{chapter6_eq32}
    \end{align}
    \begin{align}
    2V_
    {0,j}^*=(1-\dfrac{\Delta  t}{2}\{ (1-r_2)\delta_{\alpha_2 y}^{V^-} +r_2\delta_{\alpha_2 y}^{V^+}\}) V_{0,j}^{n+1}+(1+\dfrac{\Delta  t}{2}\{ (1-r_2)\delta_{\alpha_2 y}^{V^-} +r_2\delta_{\alpha_2 y}^{V^+}\}) V_{0,j}^{n}\nonumber \\
    2V_
    {N_x,j}^*=(1-\dfrac{\Delta  t}{2}\{ (1-r_2)\delta_{\alpha_2 y}^{V^-} +r_2\delta_{\alpha_2 y}^{V^+}\}) V_{N_x,j}^{n+1}+(1+\dfrac{\Delta  t}{2}\{ (1-r_2)\delta_{\alpha_2 y}^{V^-} +r_2\delta_{\alpha_{2} y}^{V^+}\}) V_{N_x,j}^{n}.
    \label{chapter6_eq32a}
    \end{align}
    Dirichlet boundary conditions are used and hence, $S_
   {0,j}^*,  S_
   {N_x,j}^*,  I_
   {0,j}^*$, $I_
   {N_x,j}^*$, $V_{0,j}^*$ and $V_
   {N_x,j}^*$  becomes zero.     
     The numerical solutions acquired by using the ADI-CN scheme for the fractional diffusion $SIV$-model is given by (\ref{chapter6_eq35C}) are shown below. To see the diffusion of the infection the following set of initial conditions are used,  
 \begin{align}
 I(0)=\left\{
    \begin{array}{c l}    
        0.1& \text{mid point of the finite grid}\\
        0 & \text{elsewhere}
    \end{array}\right. 
 \end{align}
  \begin{align}
      S(0)=\left\{
    \begin{array}{c l}    
        0.9& \text{mid point of the finite grid}\\
        1 & \text{elsewhere}
    \end{array}\right.
 \end{align}
 \begin{align}
  V(0)=\left\{
    \begin{array}{c l}    
        0 & \text{everywhere}
    \end{array}\right.
 \end{align}
 
  \noindent  In this paper Dirichlet  boundary conditions are of interest.   Dirichlet conditions for both the $S$ and $I$ compartments  on the rectangular region $x_L\leq x \leq x_H$ and $y_L \leq y \leq y_H$ are of the form:
    \begin{align}
    & S(0,y,t)= 
         S(1,y,t)=  
    S(x,0,t)= 
        S(x,1,t)= 0\nonumber \\
    &I(0,y,t)=
        I(1,y,t) =
         I(x,0,t)= 
         I(x,1,t)= 0 \\
    &V(0,y,t)= 
        V(1,y,t)=  
         V(x,0,t)= 
        V(x,1,t)= 0 \nonumber.
    \end{align}    \label{chapter6_eq34}

\begin{figure}[H]    
    \centering
    \begin{subfigure}[t]{0.24\textwidth}
        \centering
        \includegraphics[width=1\textwidth]{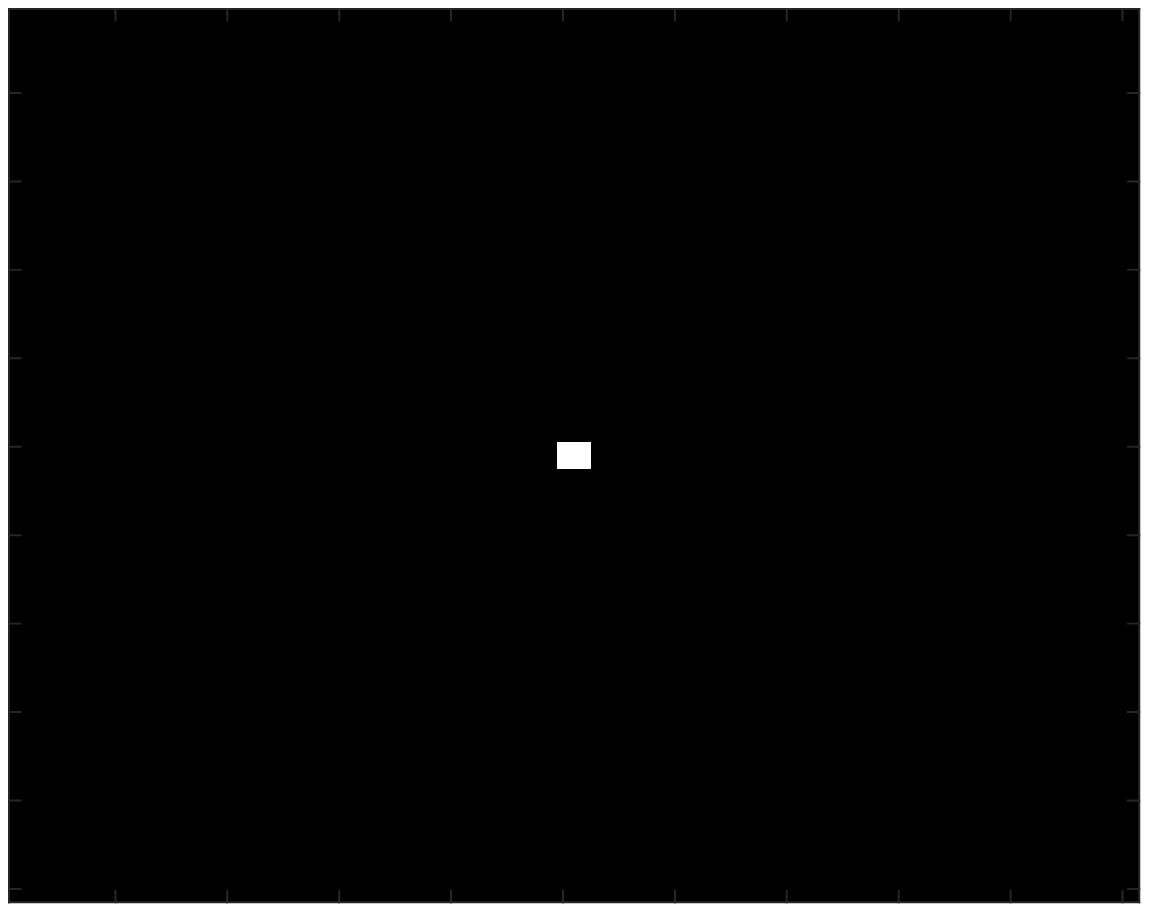}
        \caption{t= 0 days}\label{figa1}        
    \end{subfigure}
    \begin{subfigure}[t]{0.24\textwidth}
        \centering
        \includegraphics[width=1\textwidth]{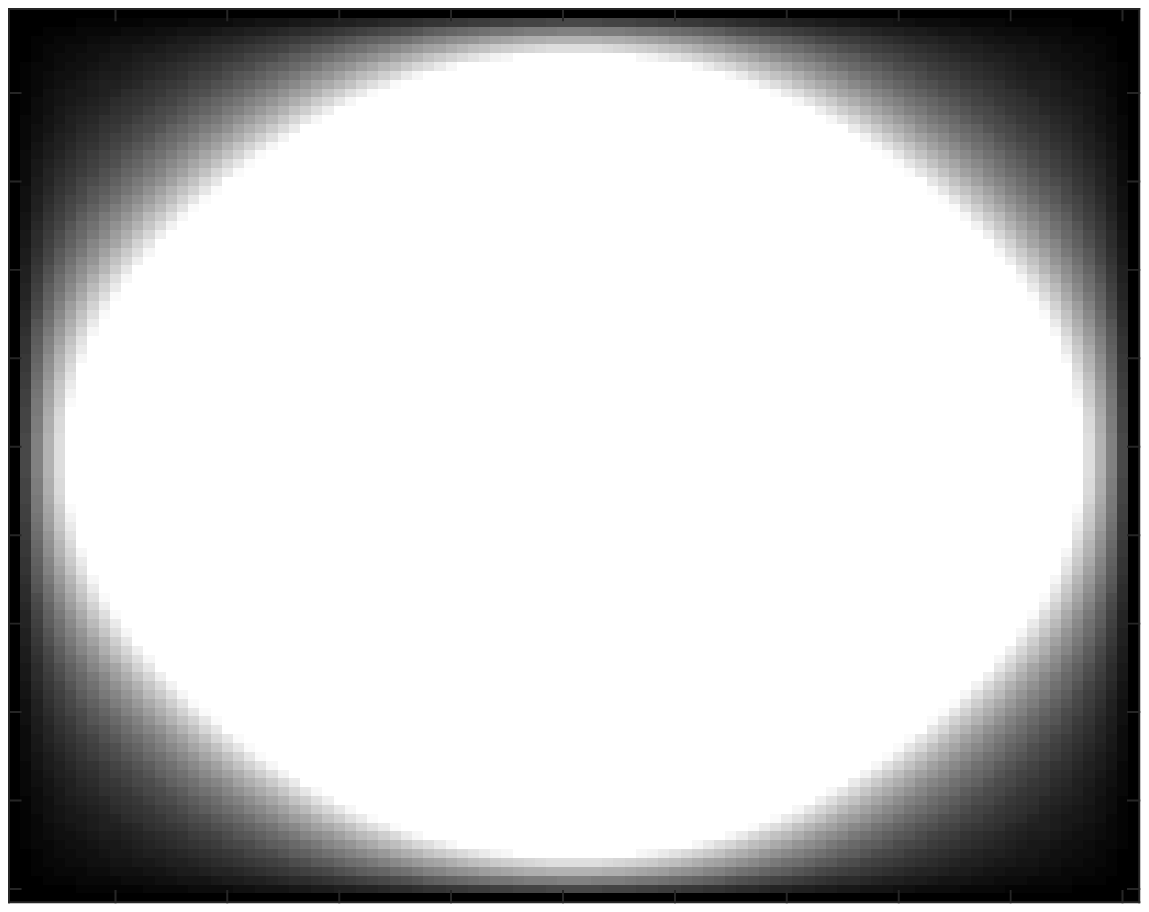}
        \caption{t= 60 days}\label{figa2}
    \end{subfigure}
    \begin{subfigure}[t]{0.24\textwidth}
        \centering
        \includegraphics[width=1\textwidth]{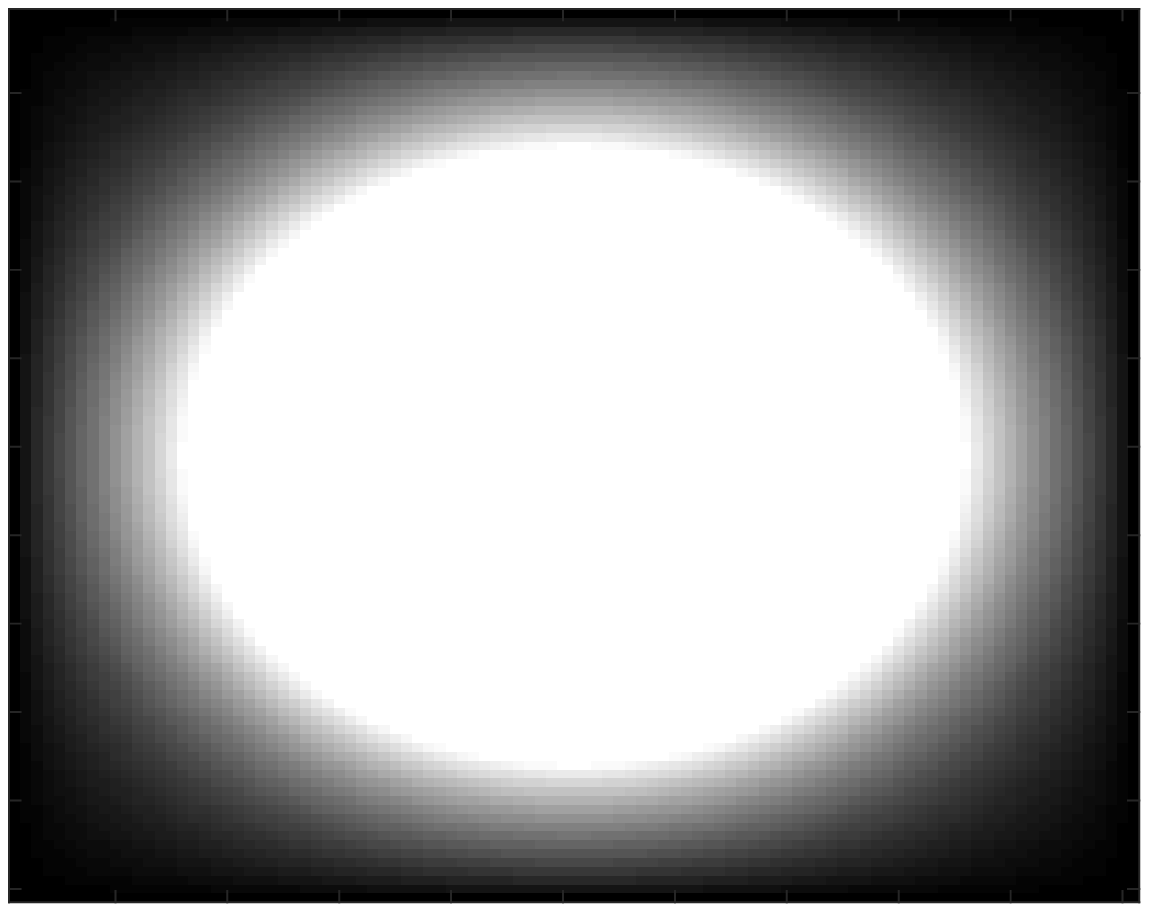}
        \caption{t= 120 days}\label{figa3}
    \end{subfigure}
    \begin{subfigure}[t]{0.24\textwidth}
    \centering
    \includegraphics[width=1\textwidth]{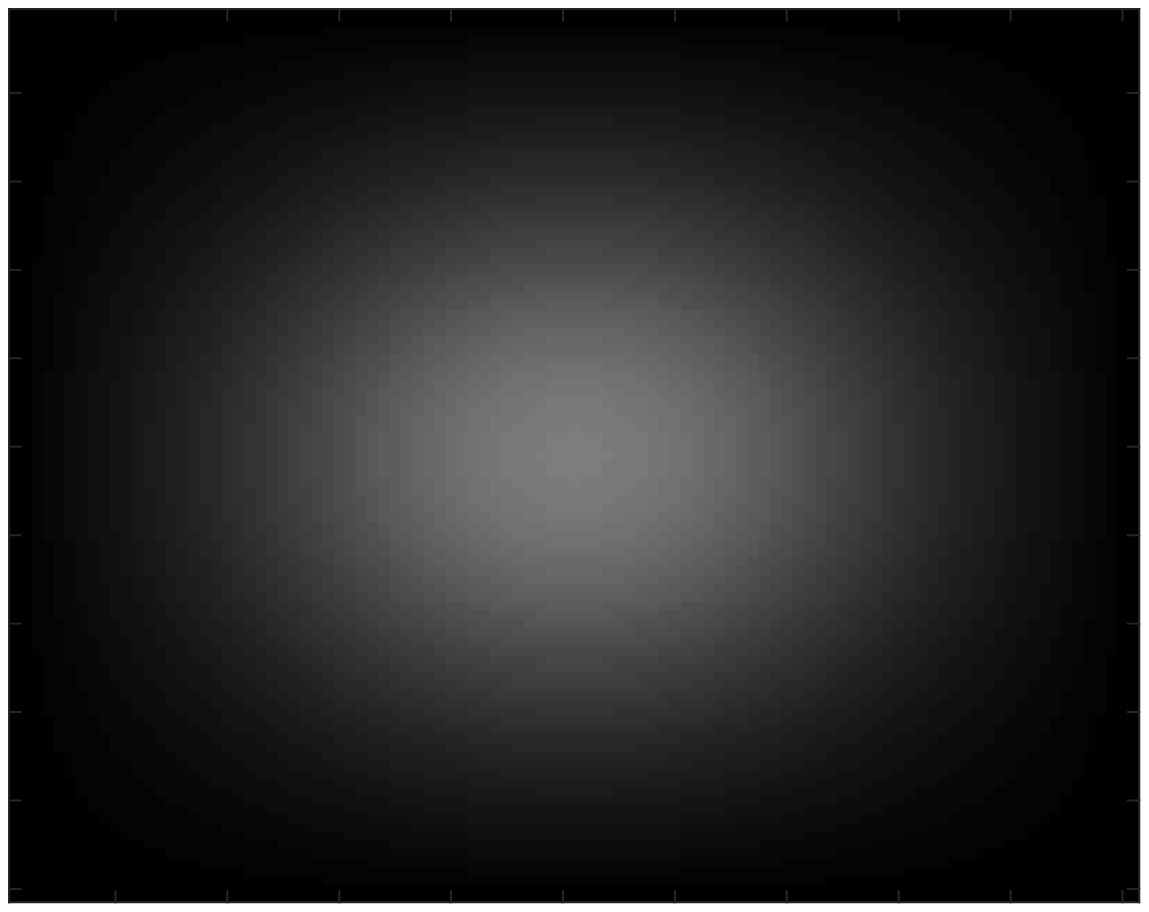}
    \caption{t= 180 days}\label{figa4}
\end{subfigure}

\begin{subfigure}[t]{0.24\textwidth}
    \centering
    \includegraphics[width=1\textwidth]{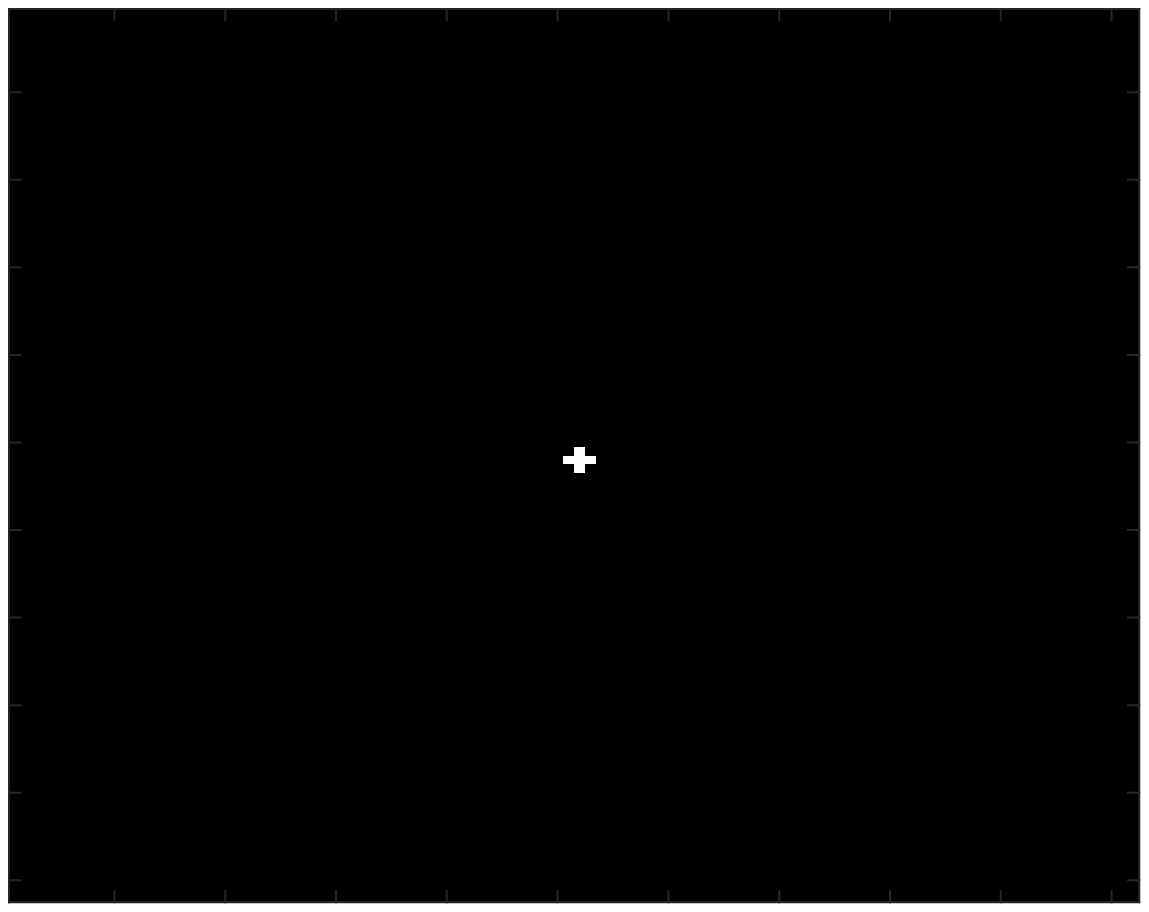}
    \caption{t= 0 days}\label{figa9}        
\end{subfigure}
\begin{subfigure}[t]{0.24\textwidth}
    \centering
    \includegraphics[width=1\textwidth]{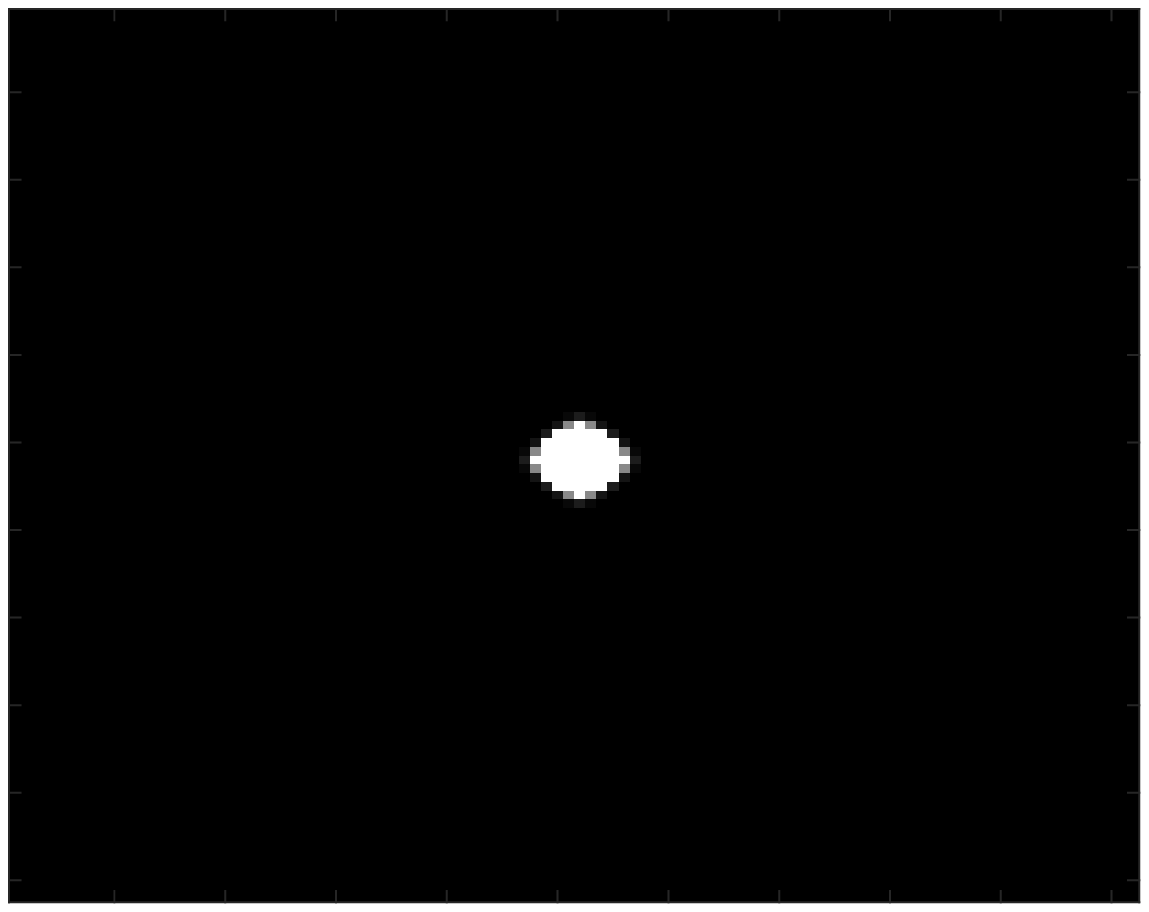}
    \caption{t= 60 days}\label{figa10}
\end{subfigure}
\begin{subfigure}[t]{0.24\textwidth}
    \centering
    \includegraphics[width=1\textwidth]{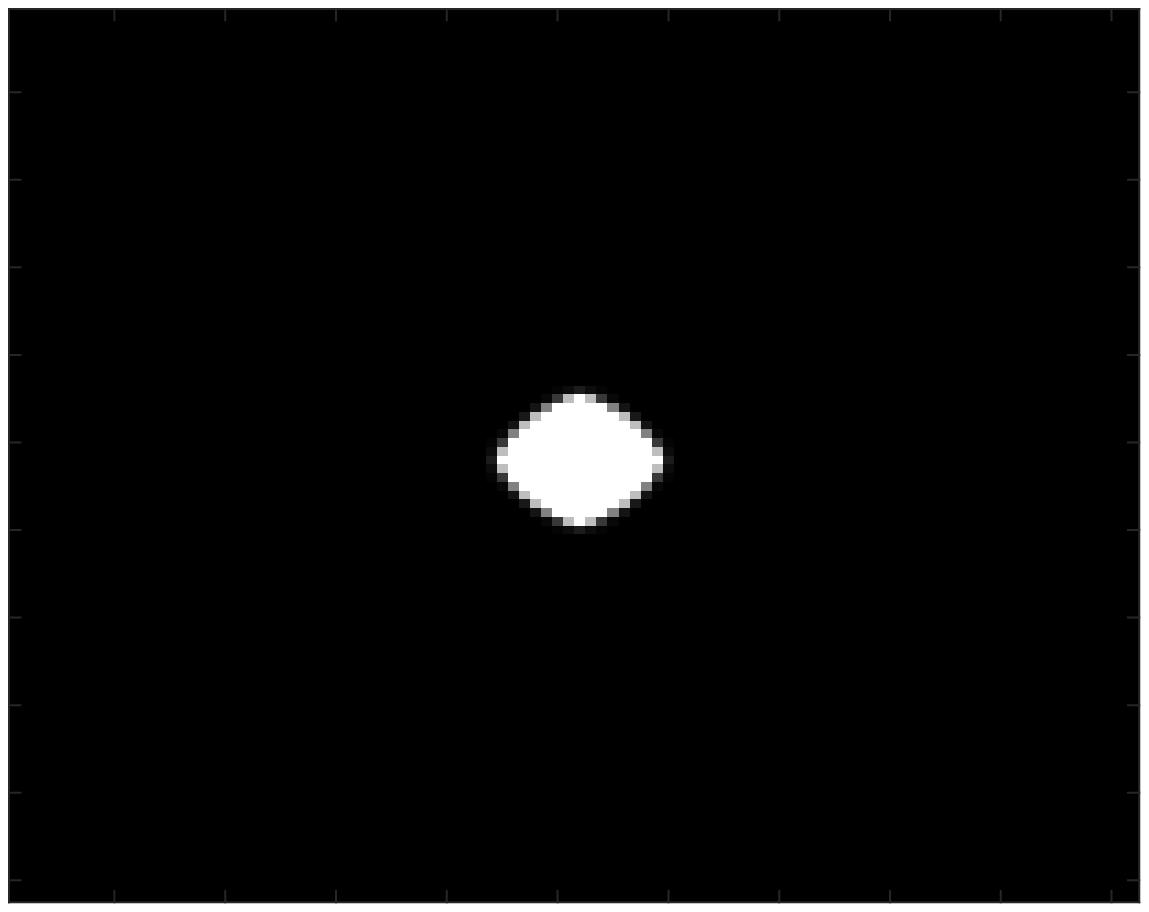}
    \caption{t= 120 days}\label{figa11}
\end{subfigure}
\begin{subfigure}[t]{0.24\textwidth}
    \centering
    \includegraphics[width=1\textwidth]{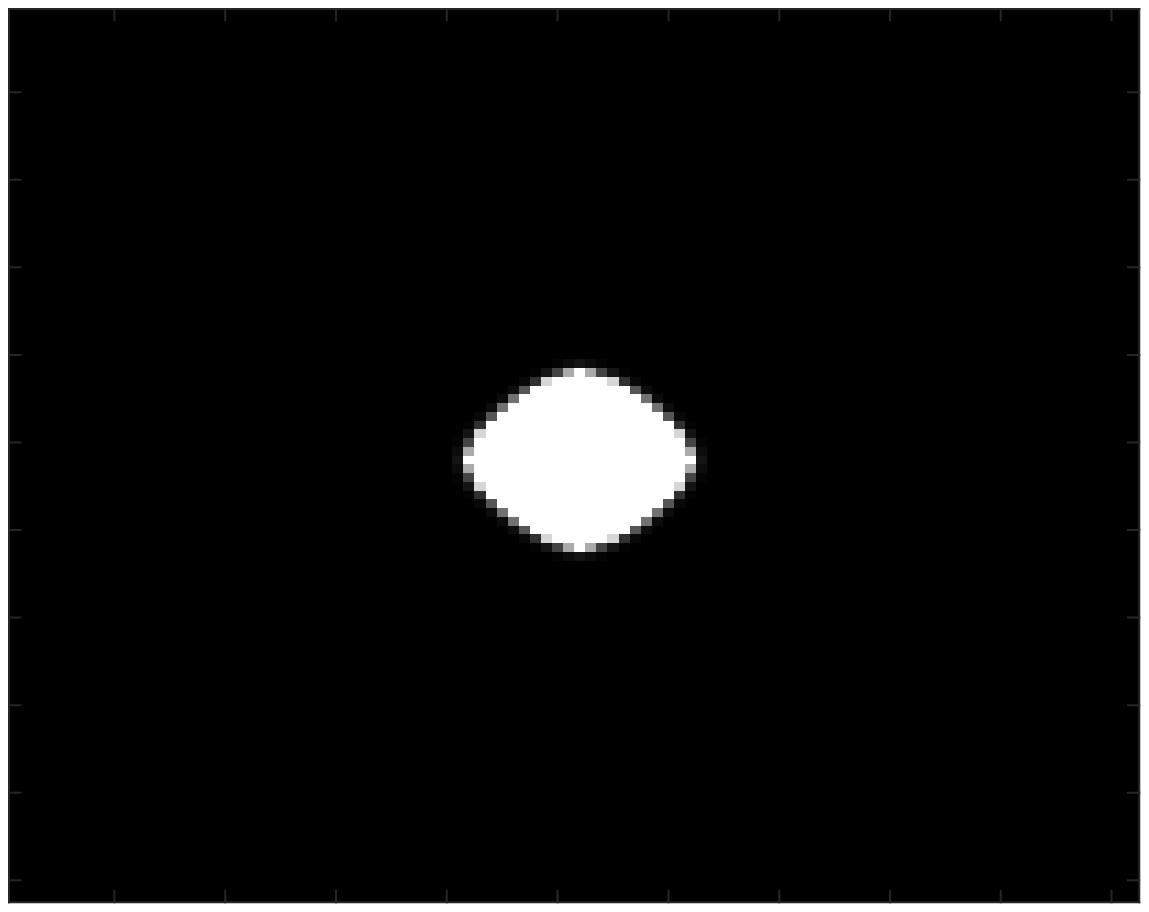}
    \caption{t= 180 days}\label{figa12}
\end{subfigure}
   \caption{Spread of the  infected hosts $I$  by Fractional Diffusion $SIV$-model where $\alpha=1.2.$ and the classical model in the host compartment.  }\label{chapter6_figure2}
\end{figure}

\begin{figure}[H]    
    \centering

    \begin{subfigure}[t]{0.24\textwidth}
        \centering
        \includegraphics[width=1\textwidth]{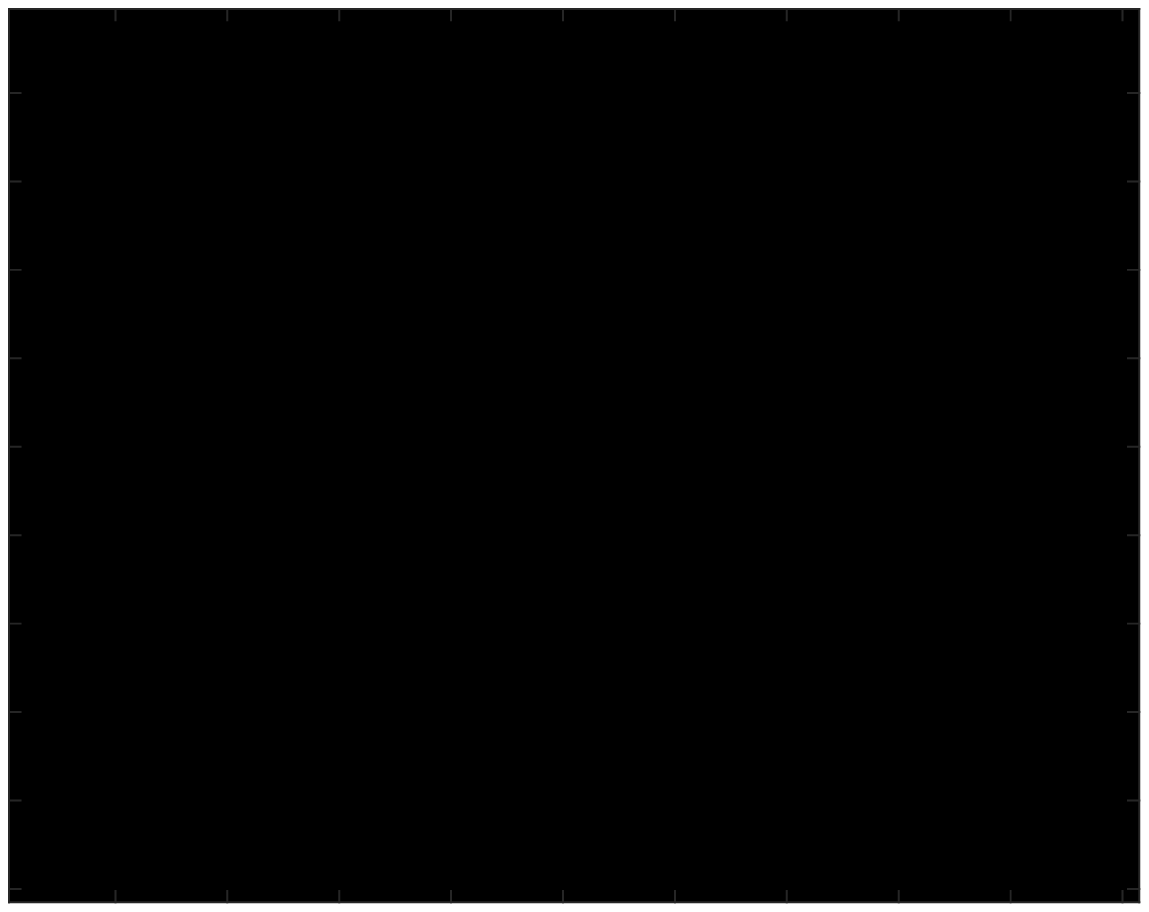}
        \caption{t= 0 days}\label{figa5}        
    \end{subfigure}
    \begin{subfigure}[t]{0.24\textwidth}
        \centering
        \includegraphics[width=1\textwidth]{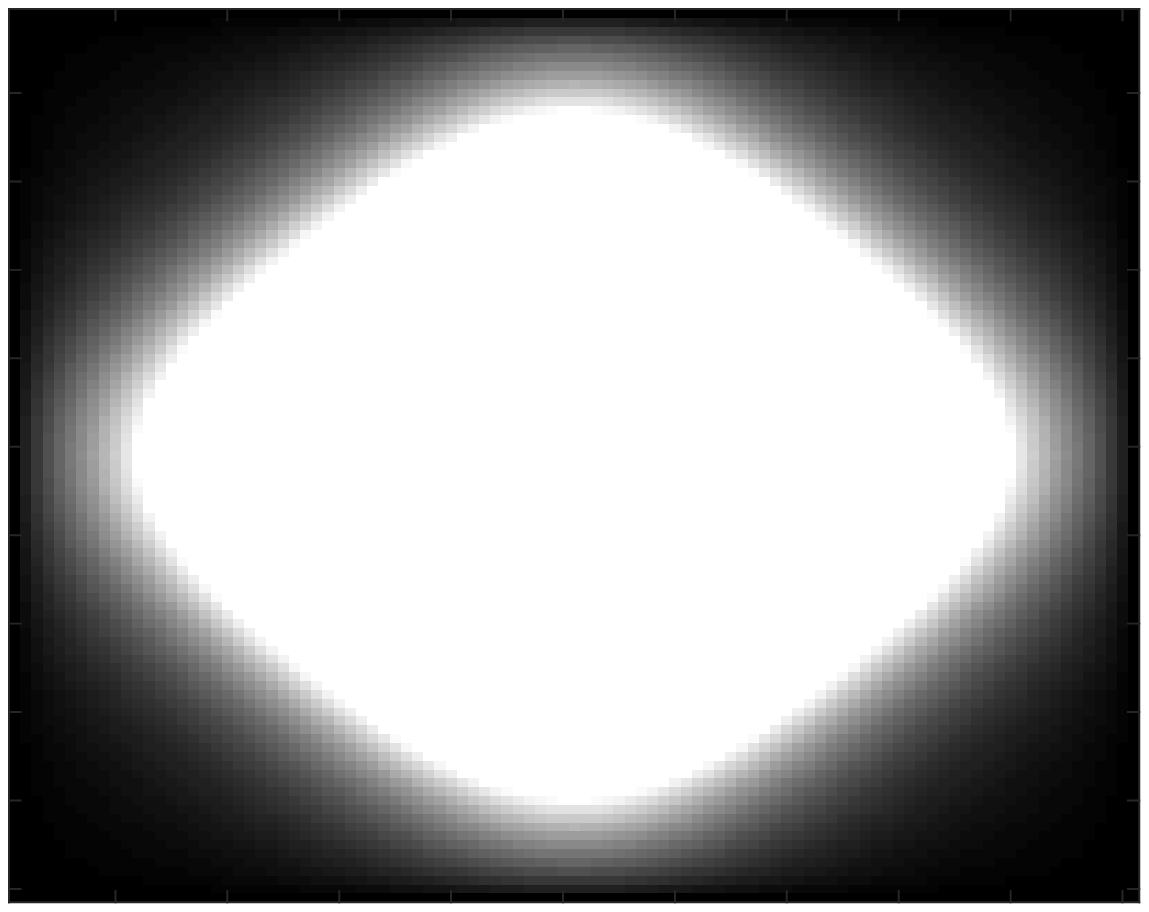}
        \caption{t= 60 days}\label{figa6}
    \end{subfigure}
\begin{subfigure}[t]{0.24\textwidth}
        \centering
        \includegraphics[width=1\textwidth]{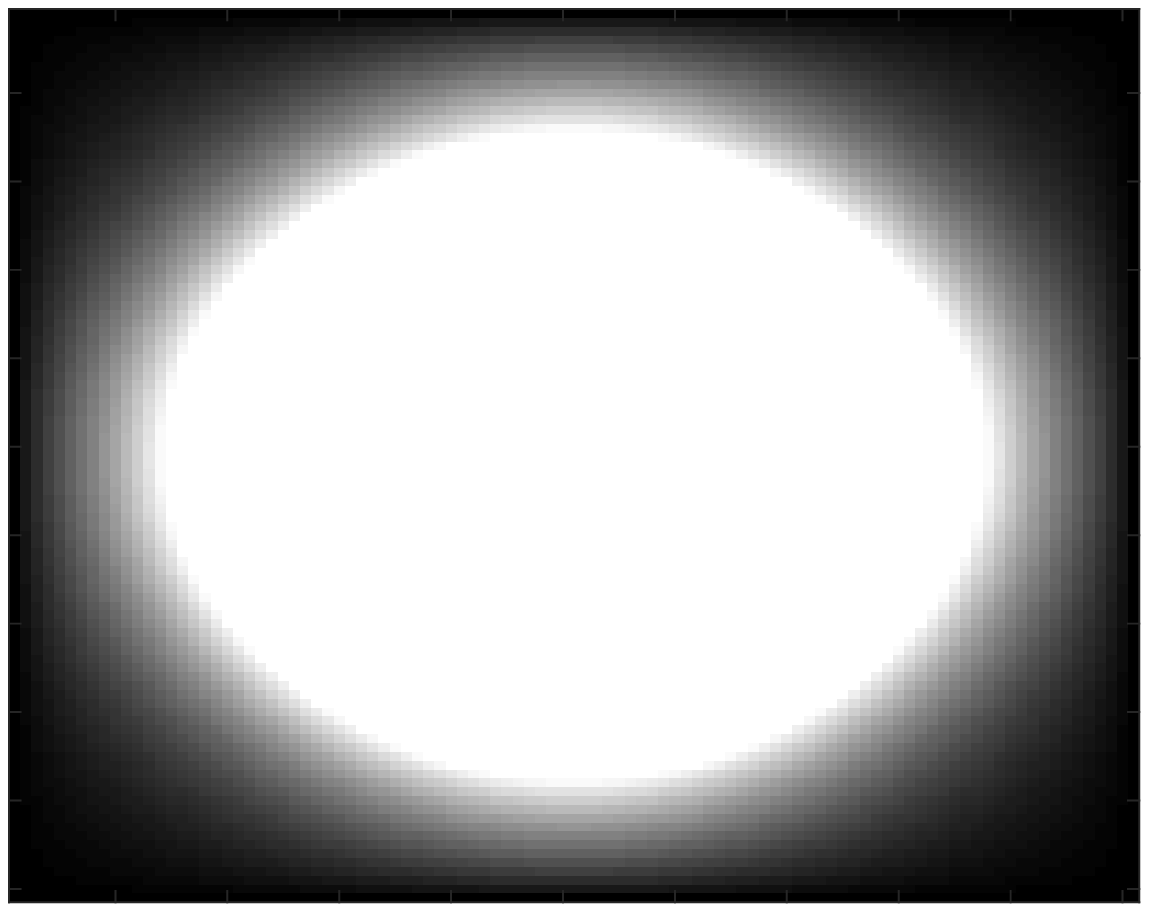}
        \caption{t= 120 days}\label{figa7}
    \end{subfigure}
    \begin{subfigure}[t]{0.24\textwidth}
        \centering
        \includegraphics[width=1\textwidth]{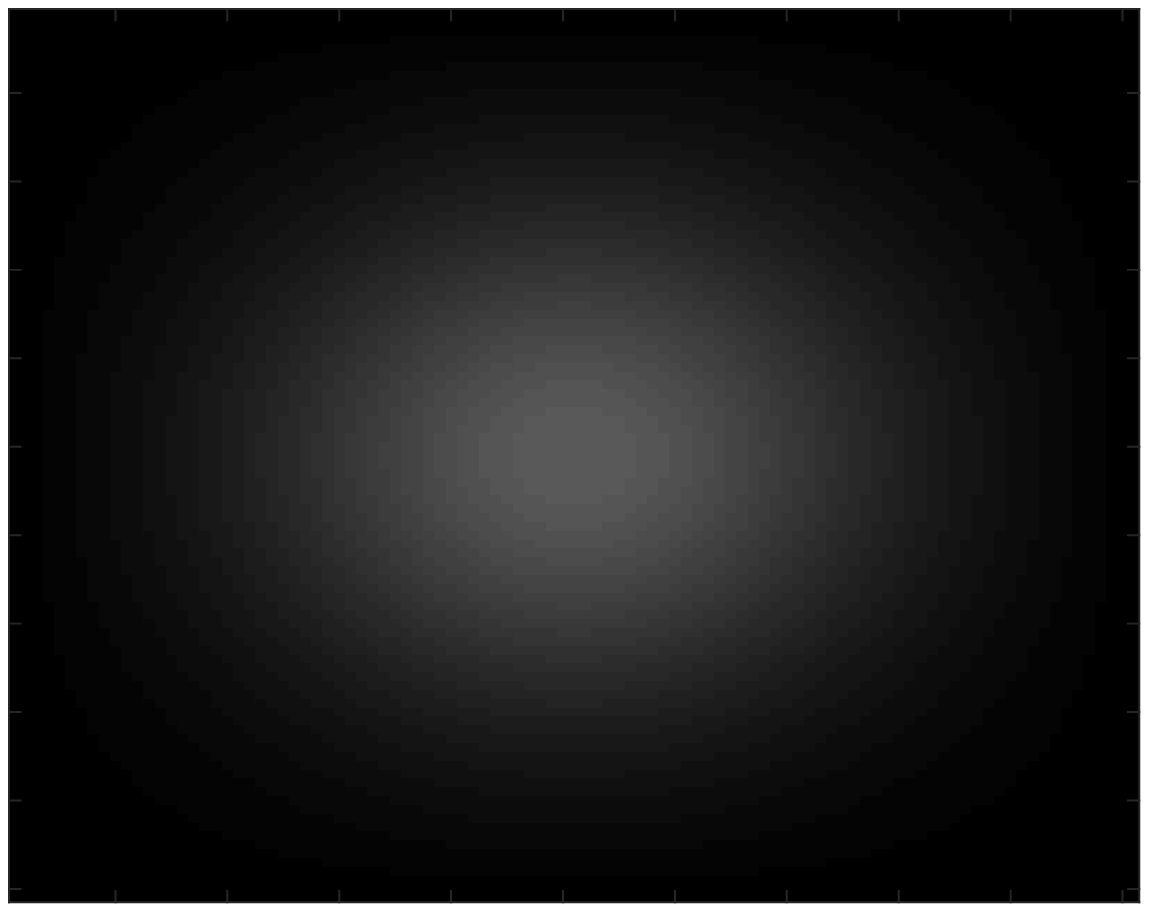}
        \caption{t= 180 days}\label{figa8}
    \end{subfigure}
    
    \begin{subfigure}[t]{0.24\textwidth}
        \centering
        \includegraphics[width=1\textwidth]{figa5}
        \caption{t= 0 days}\label{figa13}        
    \end{subfigure}
    \begin{subfigure}[t]{0.24\textwidth}
        \centering
        \includegraphics[width=1\textwidth]{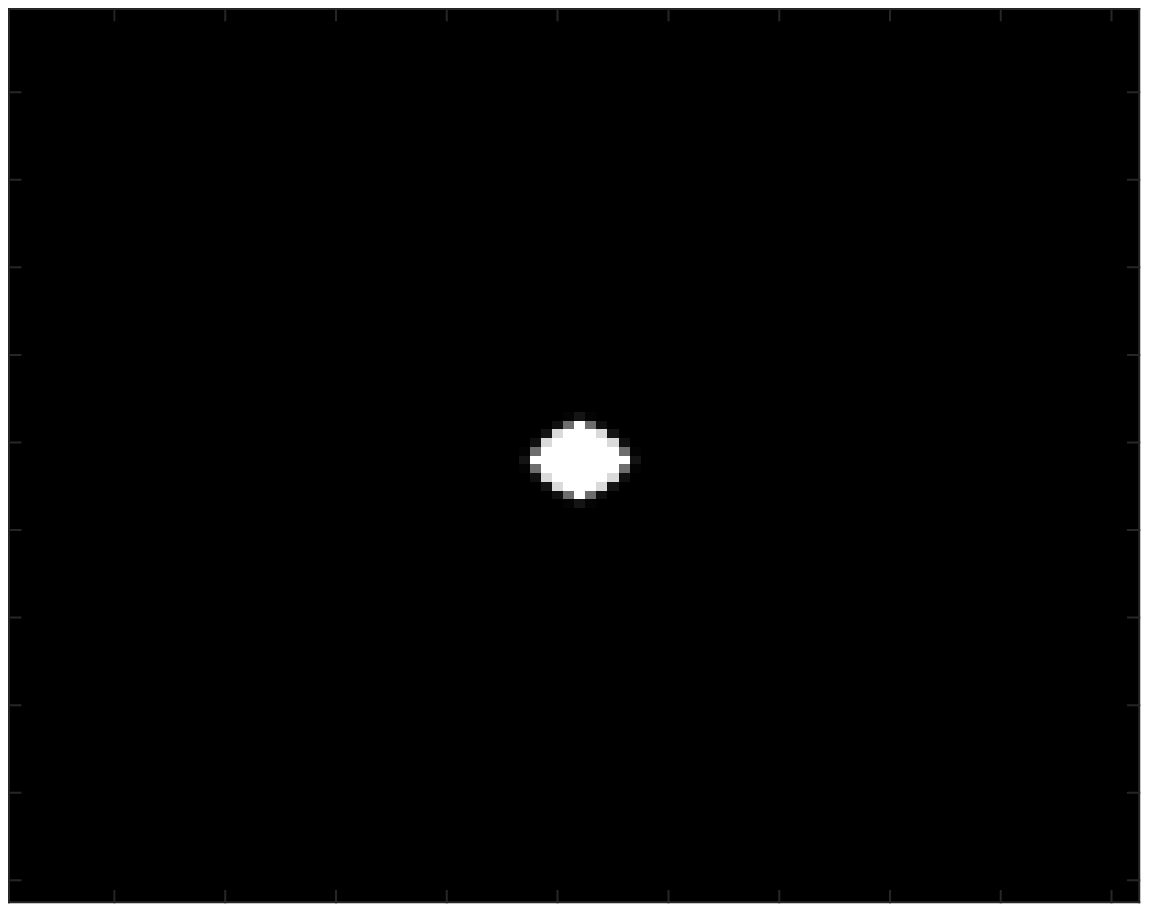}
        \caption{t= 60 days}\label{figa14}
    \end{subfigure}
 \begin{subfigure}[t]{0.24\textwidth}
        \centering
        \includegraphics[width=1\textwidth]{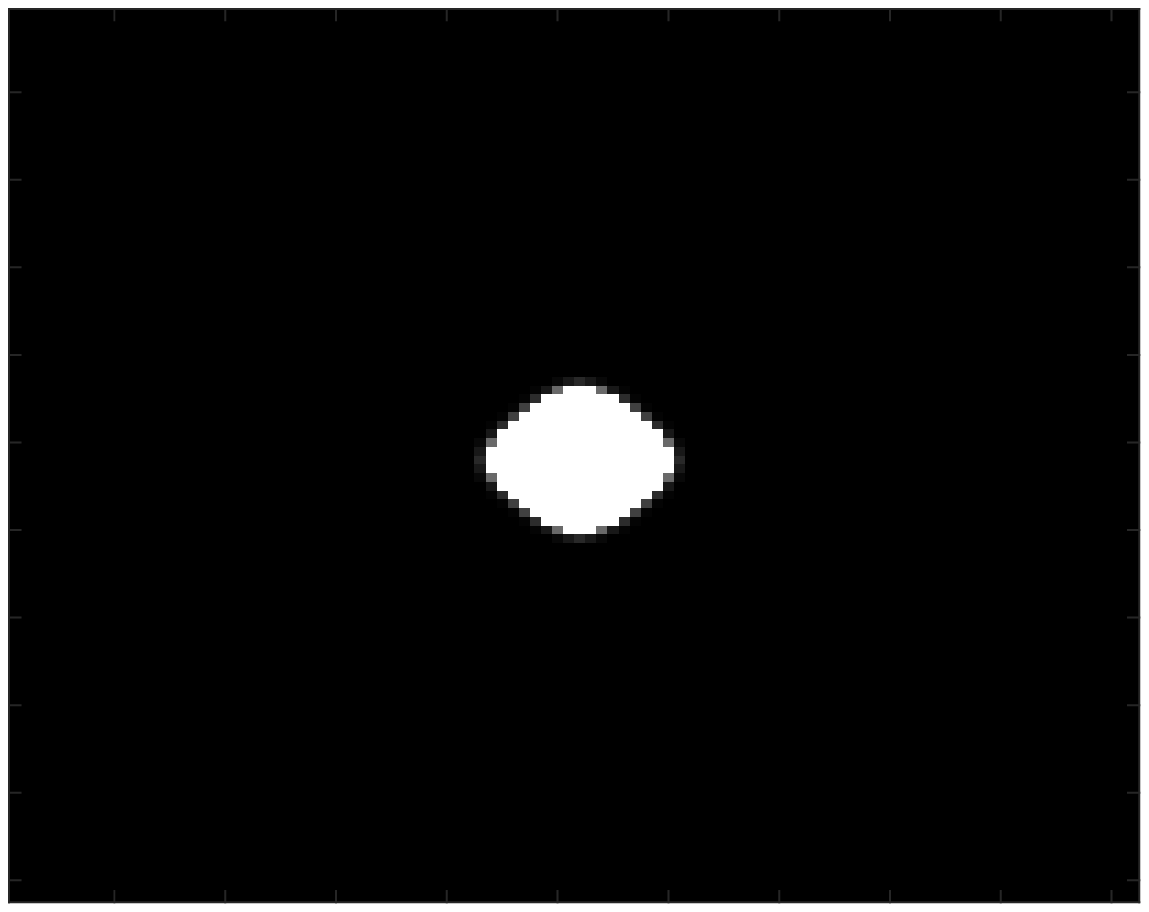}
        \caption{t= 120 days}\label{figa15}
    \end{subfigure}
    \begin{subfigure}[t]{0.24\textwidth}
        \centering
        \includegraphics[width=1\textwidth]{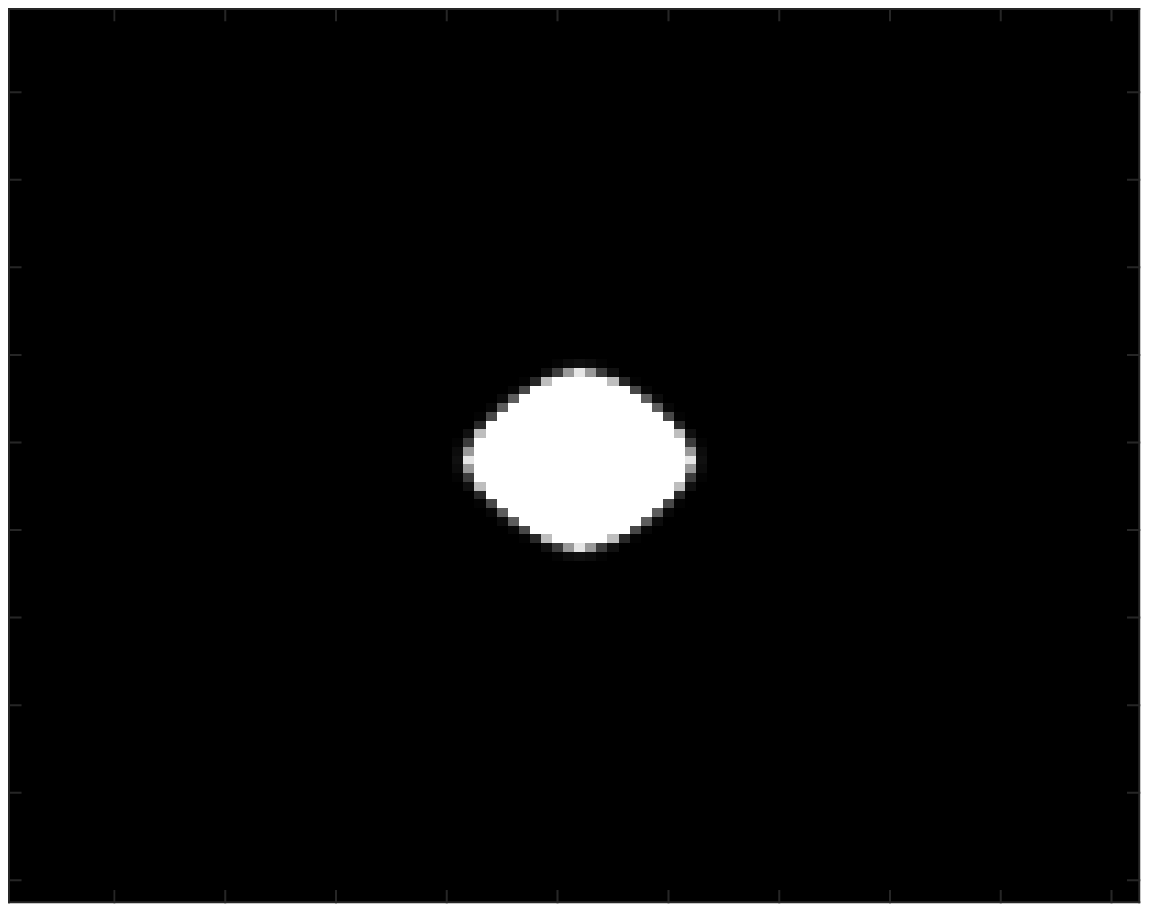}
        \caption{t= 180 days}\label{figa16}
    \end{subfigure}
    \caption{Spread of the  infected vectors $V$  by Fractional Diffusion $SIV$-model where $\alpha=1.2.$ and the classical model in the host compartment.}\label{chapter6_figure2a}
\end{figure}

    \noindent The corresponding numerical solutions of the fractional diffusion  $SIV$-model is compared with the classical diffusion $SIV$-model. Figure  \ref{chapter6_figure2} and  Figure  \ref{chapter6_figure2a} shows the results of the hosts and vectors corresponding to the fractional-order 1.2 compared with the corresponding classical case.
    
It can be seen that the spread of the classical diffusion is slower than that of the fractional diffusion.
Hence, numerical results corresponding to the fractional model shows an anomalous diffusion which can be seen in the infected hosts.

\section{Conclusion}

In this article, the reaction-diffusion $SIV$ partial differential equation model is derived by using the multi-patch system with the long term movements of the individuals. Further, here we introduce a model corresponding to the reaction-diffusion approach to the
existing $SIV$-epidemic model. The second derivative of the consistent classical
reaction-diffusion equations is substituted by using the $\alpha$ order fractional derivatives in the respective space derivatives. The model is simulated using the alternating
directions implicit (ADI) scheme with a Crank-Nicholson discretization. The
numerical results are compared with the results attained by the classical reaction-diffusion system.
The results illustrate a anomalous diffusive behaviour compared to the classical diffusion approach.

\section*{Acknowledgements}
We want to thank Thomas G\"otz and Nico Stollenwerk for intensive discussions. Moreover YJ acknowledges the support within the DAAD funded Ph.D.-scholarship program MiC at TU Kaiserslautern.

\section{Bibliography}

\end{document}